\documentclass[12pt]{article}

\usepackage{amsmath}
\usepackage{float}
\usepackage{amsfonts}
\usepackage{amstext}
\usepackage{epsf}
\usepackage[dvips]{epsfig}
\usepackage{fullpage}

\newtheorem{theorem}{Theorem}
\newtheorem{lemma}[theorem]{Lemma}

\newtheorem{fact}[theorem]{Fact}

\def \prend{\vrule depth-1pt height7pt width6pt}

\def \endpf{{\ \ \prend \medbreak}}
\def \En{{\mathbb N}}

\begin{document}

\title{Polynomial versus Exponential Growth in Repetition-Free Binary Words}

\author{
Juhani Karhum\"aki\footnote{Supported by the Academy of Finland under grant 44087.}\\
Department of Mathematics and TUCS\\
University of Turku\\
20014 Turku\\
FINLAND\\
{\tt juhani.karhumaki@cs.utu.fi}
\and
Jeffrey Shallit\footnote{Supported by a grant from NSERC.}\\
School of Computer Science \\
University of Waterloo \\
Waterloo, ON, N2L 3G1 \\
CANADA \\
{\tt shallit@graceland.uwaterloo.ca}\\
}

\date{\today}
\maketitle

\begin{abstract}
It is known that the number of overlap-free binary words of length $n$
grows polynomially, while the number of cubefree binary words grows
exponentially.  We show that the dividing line between polynomial
and exponential growth is ${7 \over 3}$.  More precisely,
there are only polynomially many binary words of length $n$ 
that avoid ${7 \over 3}$-powers, but there are exponentially many
binary words of length $n$ that avoid ${7 \over 3}^+$-powers.
This answers an open question of Kobayashi from 1986.
\end{abstract}

\section{Introduction}

We are concerned in this paper with problems on combinatorics of
words \cite{Lothaire:1983,Lothaire:2002}.

Let $\Sigma$ be a finite nonempty set, called an {\sl alphabet}.  We
consider finite and infinite words over $\Sigma$.  The set of all
finite words over $\Sigma$ is denoted by $\Sigma^*$. The set of all infinite
words (that is, maps from $\En$ to $\Sigma$) is denoted by
$\Sigma^\omega$.    In this paper we often use a particular class
of alphabets, namely
$$ \Sigma_k := \lbrace 0, 1, \ldots, k-1 \rbrace .$$

     A {\sl morphism} is a map $h: \Sigma^* \rightarrow \Delta^*$
such that $h(xy) = h(x)h(y)$ for all $x, y \in \Sigma^*$.  A morphism
may be specified by providing the image words $h(a)$ for all $a \in \Sigma$.
If $h:\Sigma^* \rightarrow \Sigma^*$ and
$h(a) = ax$ for some letter $a \in \Sigma$, then we say that
$h$ is {\sl prolongable} on $a$, and we can then iterate $h$ infinitely
often to get the fixed point
$h^\omega(a) := a \, x \, h(x) \, h^2(x) \, h^3(x) \cdots $.
If there exists an integer $k$ such that
the morphism $h$ satisfies $h(a) = k$ for all $a \in \Sigma$, we say
it is {\sl $k$-uniform}.  If a morphism is $k$-uniform for some $k$, then
we say it is {\sl uniform}.  For more on morphisms, see, for example, 
\cite{Choffrut&Karhumaki:1997}.

     A {\sl square} is a nonempty word of the form $xx$, as in the
English word {\tt murmur}.  A {\sl cube} is a nonempty word of the
form $xxx$, as in the Finnish word {\tt kokoko}.  An
{\sl overlap} is a word of the form $axaxa$, where $x$ is a possibly
empty word and $a$ is a single letter, as in the English word
{\tt alfalfa}.  A word $v$ is a {\sl factor} (sometimes called a {\sl subword})
of a word $x$ if $x$ can be written $x = uvw$ for some words $u, w$.
A word {\sl avoids squares} (resp., cubes, overlaps)
if it contains no factor 
that is a square (resp., cube, overlap).  Such words are also
called {\sl squarefree} (resp., {\sl cubefree}, {\sl overlap-free}).
For example, the English word ${\tt square}$ is squarefree, whereas
${\tt squarefree}$ is not.

      It is well-known and easily proved that every word of length
$4$ or more over a two-letter alphabet contains a square as a factor.
However, Thue proved in 1906 \cite{Thue:1906} that there exist
infinite squarefree words over a three-letter alphabet.
Thue also proved that the word $\mu^\omega (0) = 0110100110010110\cdots$
is overlap-free (and hence cubefree); here
$\mu$ is the {\tt Thue-Morse morphism}
sending $0 \rightarrow 01$ and $1 \rightarrow 10$.

     Dejean \cite{Dejean:1972} 
initiated the study of fractional
powers.  Let $\alpha$ be a rational number $ \geq 1$.  Following Brandenburg 
\cite{Brandenburg:1983}, we say that a
nonempty word $w$ is an {\sl $\alpha$-power}
if there exist words $y, y'\in \Sigma^*$
such that $w = y^n y'$, and $y'$ is a prefix of $y$ with
$n + |y'|/|y| = \alpha$.  For example, the English word
{\tt alfalfa} is a ${7 \over 3}$-power and the word
{\tt ionization} is a ${{10} \over 7}$-power.  If $\alpha$ is a real
number, we say that a word {\sl $w$ avoids $\alpha$-powers} (or is
{\sl $\alpha$-power-free}) if it contains no factor that is a $\beta$-power
for any rational $\beta \geq \alpha$.  We say that a word $w$ avoids $\alpha^+$-powers
(or is $\alpha^+$-power-free) if it contains no factor that is
a $\beta$-power for rational $\beta > \alpha$.  Thus a word is overlap-free iff
it is $2^+$-power-free.

     We may enumerate the number of words avoiding various patterns.
Brandenburg \cite{Brandenburg:1983} proved (among other things) that there are exponentially
many cubefree binary words; also see Edlin \cite{Edlin:1999}.
Restivo and Salemi \cite{Restivo&Salemi:1983,Restivo&Salemi:1985a} proved
that there exist only polynomially many
overlap-free binary words of length
$n$; in fact they gave an upper bound of $O(n^{\log_2 15})$.
The exponent $\log_2 15$ was improved to $1.7$ by Kfoury \cite{Kfoury:1988a}
to $1.587$ by Kobayashi \cite{Kobayashi:1988}, and to $1.37$ by 
Lepist\"o \cite{Lepisto:1995}.  Also see Cassaigne 
\cite{Cassaigne:1993b} 

     Overlap-free words avoid $2^+$-powers, and there are only
polynomially many over $\Sigma_2$.
Cubefree words avoid $3$-powers and there are exponentially many
over $\Sigma_2$.
Kobayashi \cite[Problem 6.6]{Kobayashi:1986} asked the following
natural question: at what exponent $\alpha$ (if any) does 
the number of binary words avoiding $\alpha$-powers jump from polynomial
to exponential?  In this paper we prove that the answer is $7 \over 3$.
Our proof uses the fact that various structure theorems, which previously were
known for overlap-free words, also hold for any exponent $\alpha$
with $2 < \alpha \leq {7 \over 3}$.

\section{Preliminary lemmas}

     We begin with some notation and preliminary results.
We write $\overline{0} = 1$ and $\overline{1} = 0$.  We let
$\mu$ be the Thue-Morse morphism mentioned in the previous
section.

\begin{lemma}
Let $t, v \in \Sigma_2^*$.
If there exist $c, d \in \Sigma_2$
such that $c \mu(t) = \mu(v) d$, then $d = c$,
$t = \overline{c}^n$, and $v = c^n$, where $n = |t| = |v|$.
\label{as1}
\end{lemma}

\begin{proof}
See \cite[Lemma 1.7.2]{Allouche&Shallit:2003}.
\endpf
\end{proof}

\begin{lemma}  
Suppose $t, y \in \Sigma_2^*$ and $\mu(t) = yy$.  Then there exists
$v \in \Sigma_2^*$ such that $y = \mu(v)$.
\label{as2}
\end{lemma}

\begin{proof}
See \cite[Lemma 1.7.3]{Allouche&Shallit:2003}.
\endpf
\end{proof}

\begin{lemma}
Let $h:\Sigma^* \rightarrow \Sigma^*$ be a uniform morphism, and let
$\alpha$ be a rational number.  If $w$ contains an $\alpha$-power,
then $h(w)$ contains an $\alpha$-power.
\label{as5}
\end{lemma}

\begin{proof}
     Suppose $w$ contains an $\alpha$-power.  Then there exist
words $s, s' \in \Sigma^+$ and $r, t \in \Sigma^*$ such that
$w = r s^n s' t$, where $s'$ is a nonempty prefix of $s$ and
$n + |s'|/|s| = \alpha$.
Then $h(w) = h(r) h(s)^n h(s') h(t)$.  Then $h(w)$ 
contains the $\alpha$-power $h(s)^n h(s')$.
\endpf
\end{proof}

     Note that for arbitrary morphisms the result need not be true (unless
$\alpha$ is an integer).

\begin{lemma}
      Let $w \in \Sigma_2^*$, and suppose $\mu(w)$ contains
an $\alpha$-power.  Then $w$ contains a $\beta$-power with
$\beta \geq \alpha$.
\label{as6}
\end{lemma}

\begin{proof}
     Suppose $\mu(w)$ contains an $\alpha$-power, say
$\mu(w) = x y^n y' z$, where $n + |y'|/|y| = \alpha$.  
There are four cases to consider, based on the parity of $|x|$ and $|y|$.

\medskip

     Case 1:  $|x|$ is even and $|y|$ is even.  There are two subcases,
     depending on the parity of $|y'|$.

     \medskip

     	Case 1a:  $|y'|$ is even. Then $|z|$ is even.  Then
     	there exist words $r, s, s', t$, with $s'$ a prefix of $s$, such that
     	$\mu(r) = x$, $\mu(s) = y$, $\mu(s') = y'$, and
     	$\mu(t) = z$.  Then $w = r s^n s' t$, and so $w$ contains
     	the $\alpha$-power $s^n s'$.

	\medskip

	Case 2a:  $|y'|$ is odd.  Then $|z|$ is odd.  Then 
	there exist words $r, s, s', t$, with $s'$ a prefix of $s$,
	and a letter $c$ such
	that $\mu(r) = x$, $\mu(s) = y$, $\mu(s')c = y'$, and
	$\overline{c} \mu(t) = z$.      Since $|y'|$ is odd, $|y|$ is even, and
	$y'$ is a prefix of $y$, it follows that $y' \overline{c}$ is also a prefix
	of $y$.  Hence $s' c$ is a prefix of $s$.   Then
	$w$ contains the $\beta$-power $s^n s' c$, where
	$$ \beta = n + {{|s' c|} \over {|s|}} = n + {{2 |s'| + 2} \over {2|s|}}
		= n + {{|y'| + 1} \over {|y|}} > n + {{|y'|} \over {|y|}} \geq \alpha.$$

\bigskip

     Case 2:  $|x|$ is even and $|y|$ is odd.  Then 
     there exists a word $t$ such that $\mu(t) = yy$.
     From Lemma~\ref{as2} there exists $v$ such that $y = \mu(v)$.
     But then $|y|$ is even, a contradiction. Thus this case cannot occur.

\bigskip

     Case 3:  $|x|$ is odd and $|y|$ is even.  There are two subcases,
     depending on the parity of $|y'|$.

     \medskip

     	Case 3a:  $|y'|$ is even.  Then $|z|$ is odd.   Then there exist
	words $r, s, s', t$ and letters $c, d, e$ such that
	$x = \mu(r) c$, $y = \overline{c} \mu(s) d$, $y' = \overline{d} \mu(s') e$,
	and $z = \overline{e} \mu(t)$.   Consideration of the factor $yy$
	gives $c = d$.   Hence $\mu(w) = \mu(r (cs)^ n cs' et)$
	and so $w = r (cs)^n cs' et$.  Thus $w$ contains the
	$\alpha$-power $(cs)^n cs'$, and since $s'$ is a prefix of $s$, it
	follows that $cs'$ is a prefix of $cs$.

	\medskip

	Case 3b:  $|y'|$ is odd.  Then $|z|$ is even.  
	Then we are in the mirror image of case 2a, and the same proof works.


	\medskip

     Case 4:  $|x|$ is odd and $|y|$ is odd.  Then from length considerations
     we see that there exist words $t, v$ and letters $c, d$ such that
     $y = c \mu(t) = \mu(v) d$.  By Lemma~\ref{as1}, we have
     $d = c$, $t = {\overline{c}}^n$, $v = c^n$.     Thus
     $y = c (\overline{c} c)^n$.
     Since $y'$ is a nonempty prefix of $y$, we may write
     $\mu(w) = x y^2 c t$ for some word $t$.
     Since $|x|$ and $|y|$ are odd, and $y$ ends in $c$, we must have
     that $cc$ is the image of letter under $\mu$, a contradiction.
     Thus this case cannot occur.
\endpf
\end{proof}

\begin{theorem}
Let $w \in \Sigma_2^*$, and let $\alpha > 2$ be a real number.
Then $w$ is $\alpha$-power-free iff $\mu(w)$ is $\alpha$-power-free.
\label{as3}
\end{theorem}

\begin{proof}
    Combine Lemmas~\ref{as5} and \ref{as6}.
\endpf
\end{proof}

\bigskip

      We remark that Theorem~\ref{as3} is not true if $\alpha = 2$;
for example $w = 01$ contains no square, but $\mu(w) = 0110$ does.

\section{A structure theorem for $\alpha$-power-free words for
$2 < \alpha \leq {7 \over 3}$}

      Restivo and Salemi \cite{Restivo&Salemi:1983,Restivo&Salemi:1985a}
proved a beautiful structure theorem for overlap-free binary words.  Roughly speaking,
it says that any overlap-free word is, up to removal of a short prefix or suffix,
the image of another overlap-free word under $\mu$, the Thue-Morse morphism.
Perhaps surprisingly,
the same sort of structure theorem exists for binary words avoiding
$\alpha$-powers, where $\alpha$ is any real number with
$2 < \alpha \leq {7 \over 3}$.

\begin{theorem}
      Let $x$ be a word avoiding $\alpha$-powers, with $2 < \alpha \leq {7 \over 3}$.
Let $\mu$ be the Thue-Morse morphism.
Then there exist $u, v, y$ with $u, v \in \lbrace \epsilon, 0, 1, 00, 11 \rbrace$
and a word $y \in \Sigma_2^*$ avoiding $\alpha$-powers, such that $x = u \mu(y) v$.
\label{rs}
\end{theorem}

\begin{proof}
    We prove the result by induction on $|x|$.  If $|x| \leq 2$, then 
the factorizations can be chosen as shown in the following table.
\begin{center}
\begin{tabular}{|r|r|r|r|}
\hline
$x$ & $u$ & $y$ & $v$ \\
\hline
$\epsilon$ & $\epsilon$ & $\epsilon$ & $\epsilon$ \\
$0$ & $0$ & $\epsilon$ & $\epsilon$ \\
$1$ & $1$ & $\epsilon$ & $\epsilon$ \\
$00$ & $00$ & $\epsilon$ & $\epsilon$ \\
$01$ & $\epsilon$ & $0$ & $\epsilon$ \\
$10$ & $\epsilon$ & $1$ & $\epsilon$ \\
$11$ & $11$ & $\epsilon$ & $\epsilon$ \\
\hline
\end{tabular}
\end{center}

      Now suppose the claim is true for all $x$ with $|x| < k$.  We prove
it for $|x| = k$.  Let $x$ be $\alpha$-power-free with $|x| \geq 3$.
Write $x = az$ with $a \in \Sigma_2$ and $z \in \Sigma_2^*$.  Since
$x$ is $\alpha$-power-free, so is $z$.  Since $|z| < |x|$, by induction
there exist $u', v' \in \lbrace \epsilon, 0, 1, 00, 11 \rbrace $ 
and a $\alpha$-power-free word $y'$ such that $z = u' \mu(y') v'$.  

     Now there are several cases to consider.

\begin{itemize}

     \item[Case] 1:  $u' = \epsilon$ or $u' = a$.  Then we may write
     $x = u \mu(y) v$ with $(u, y, v) = (au', y', v')$.

     \item[Case] 2:  $u' = \overline{a}$.  Then 
     $x = u \mu(y) v$ with $(u, y, v) = (\epsilon, ay', v')$.
     Since $x$ is $\alpha$-power-free, so is $\mu(ay')$, and hence,
     by Theorem~\ref{as3}, so is $ay'$.

     \item[Case] 3:  $u' = aa$.  Then $x$ begins with $aaa = a^3$, and so $x$ does
     not avoid $\alpha$-powers.

     \item[Case] 4:  $u' = \overline{a} \, \overline{a}$.
     Then $x = a \overline{a} \, \overline{a} \mu(y') v'$.

     \begin{itemize}

     \item[Case] 4.a:  $|y'| = 0$.  Then $x = a \overline{a} \, \overline{a} v'$.
     If $v' = \epsilon$ (resp., $v' = a$, $v' = aa$), then we can write
     $x = u \mu(y) v$ with $(u, y, v) = (\epsilon, a, \overline{a})$
     (resp., $(u,y,v) = (\epsilon, a \overline{a}, \epsilon)$,
     $(u,y,v) = (\epsilon, a \overline{a}, a)$.  Otherwise,
     if $v'= \overline{a}$ or $v' = \overline{a} \, \overline{a}$
     then $x$ contains
     $\overline{a} \, \overline{a} \, \overline{a} = {\overline{a}}^3$, and
     so $x$ does not avoid $\alpha$-powers.

%
%
%
%
%

     \item[Case] 4.b:  $|y'| \geq 1$.  There are two cases to consider.

     \begin{itemize}

     \item[Case] 4.b.i:  $y' = a y''$.  There are several cases to consider.

	\begin{itemize}

	\item[Case] 4.b.i.1:  $|y''| = 0$.  Then $y' = a$ and
	$x = a \overline{a} \, \overline{a} a \overline{a} v'$.

	If $v' = \epsilon$ (resp., $v' = a$, $v'= aa$, $v' = \overline{a}$),
	then we can write $x = u \mu(y) v$ with
	$(u,y,v) = (\epsilon, a \overline{a}, \overline{a})$
	(resp., $(u, y, v) = (\epsilon, a \overline{a} \, \overline{a}, \epsilon)$,
	$(u, y, v) = (\epsilon, a \overline{a} \, \overline{a}, a)$,
	$(u, y, v) = (\epsilon, a \overline{a}, \overline{a} \,
				 \overline{a} \, )$).
	Otherwise, if $v' = \overline{a} \, \overline{a}$,
	then $x$ contains
	$\overline{a} \, \overline{a} \, \overline{a} = {\overline{a}}^3$, and
	so $x$ does not avoid $\alpha$-powers.

%
%
%


	\item[Case] 4.b.i.2:  $|y''| \geq 1$.  
	If $y'' = a y'''$,
	then $x = a \overline{a} \,
	\overline{a} a \overline{a} a \overline{a} \mu(y''') v'$,
	so $x$ contains the ${5 \over 2}$-power
	$\overline{a} a \overline{a} a \overline{a}$.
        If $y'' = \overline{a} y'''$, then
	$x =  a \overline{a} \,
	  \overline{a} a \overline{a} \, \overline{a} a \mu(y''') v'$,
	so $x$ contains the ${7 \over 3}$-power
	$a \overline{a} \, \overline{a} a \overline{a} \, \overline{a} a$.



	\end{itemize}

	\item[Case] 4.b.ii:  $y' = \overline{a} y''$. 
     Then $x = a \overline{a} \, \overline{a} \, \overline{a} a \mu(y'') v'$.
	Thus $x$ contains
	     $\overline{a} \, \overline{a} \, \overline{a} = {\overline{a}}^3$.

	\end{itemize}

     \end{itemize}

\end{itemize}

Our proof by induction is now complete.

\endpf
\end{proof}

     The decomposition in Theorem~\ref{rs} is actually unique
if $|x| \geq 7$.  As this requires more tedious case analysis and is not
crucial to our discussion, we do not prove this here.

     We also note that the role of $7 \over 3$ in Theorem~\ref{rs} is crucial,
since no word $\cdots 0110110 \cdots$ can be factorized
in the stated form.

\section{Polynomial upper bound on the number of $7 \over 3$-power-free words}

      Theorem~\ref{rs} has the following implication.  
Let $x = x_0$ be a nonempty binary word that is
$\alpha$-power-free, with $2 < \alpha \leq {7 \over 3}$.  Then
by Theorem~\ref{rs} we can write $x_0 = u_1 \mu(x_1) v_1$ 
with $|u_1|, |v_1| \leq 2$.  If $|x_1| \geq 1$, we can repeat
the process, writing
$x_1 = u_2 \mu(x_2) v_2$.  Continuing in this fashion, we obtain
the decomposition $x_i = u_i \mu(x_i) v_i$ 
until $|x_{t+1}| = 0$ for some $t$.  Then
$$ x_0 = u_1 \mu(u_2) \cdots \mu^{t-1} (u_{t-1}) \mu^t (x_t)
	\mu^{t-1} (v_{t-1}) \cdots \mu(v_2) v_1.$$
Then from the inequalities
$1 \leq |x_t| \leq 4$ and
$2 |x_i| \leq |x_{i-1}| \leq 2|x_i| + 4$, $1 \leq i \leq t$,
an easy induction gives
$2^t \leq |x| \leq 2^{t+3} - 4$.    Thus $t \leq \log_2 |x| < t+3$, and
so 
\begin{equation}
\log_2 |x| - 3 < t \leq \log_2 |x|.
\label{ineq1}
\end{equation}

      There are at most $5$ possibilities for each $u_i$ and $v_i$, and
there are at most $22$ possibilities for $x_t$ (since $1 \leq |x_t| \leq 4$
and $x_t$ is $\alpha$-power-free).  Inequality~(\ref{ineq1}) shows there
are at most $3$ possibilities for $t$.  Letting $n = |x|$, we see
there are at most 
$3 \cdot 22 \cdot 5^{2 \log_2 n} = 66 n^{\log_2 25}$ words of length $n$
that avoid $\alpha$-powers.    We have therefore proved

\begin{theorem}
      Let $2 <\alpha \leq {7 \over 3}$.
There are $O(n^{\log_2 25}) = O(n^{4.644})$ binary words of length $n$
that avoid $\alpha$-powers.
\label{bnd}
\end{theorem}

     We have not tried to optimize the exponent in Theorem~\ref{bnd}.  Probably it can be
made significantly smaller.

\section{Exponential lower bound on the number of ${7 \over 3}^+$-power-free words}
\label{expo-sec}

     In this section we prove that there are exponentially many
binary words of length $n$ avoiding ${7\over3}^+$-powers.

     Define the $21$-uniform morphism $h:\Sigma_4^* \rightarrow \Sigma_2^*$
as follows:

\begin{eqnarray*}
h(0) &=&  011010011001001101001  \\
h(1) &=&  100101100100110010110  \\
h(2) &=&  100101100110110010110  \\
h(3) &=&  011010011011001101001  .
\end{eqnarray*}

     We first show

\begin{lemma}
     Let $w$ be any squarefree word over $\Sigma_4$.  
Then  
\begin{itemize}
\item[(i)] $h(w)$ contains no square $yy$ with $|y| > 13$; and

\item[(ii)] $h(w)$ contains no ${7\over 3}^+$-powers.
\end{itemize}
\label{seven}
\end{lemma}

\begin{proof}
     We first prove (i).  
We argue by contradiction.  Let $w = a_1 a_2 \cdots a_n$ be a
squarefree word such that $h(w)$ contains a square,
i.e., $h(w) = xyyz$ for some $x, z \in \Sigma_4^*$,
$y \in \Sigma_4^+$.
Without loss of generality, assume that $w$ is a shortest such
word, so that $0 \leq |x|, |z| < 21$.  

Case 1:  $|y| \leq 42$.    In this case we can take $|w| \leq 5$.
To verify that $h(w)$ contains no squares $yy$ with $|y| > 13$,
it therefore suffices to check the image of 
each of the 264 squarefree words in $\Sigma_4^5$.

Case 2: $|y| > 42$.  First, we observe the following facts about $h$.

\begin{fact}
\begin{itemize}
\item[(i)]
Suppose $h(ab) = t h(c) u$ for some
letters $a, b, c \in \Sigma_4$
and words $t, u \in \Sigma_2^*$.
Then this inclusion is trivial (that is,
$t = \epsilon$ or $u = \epsilon$) or $u$ is not a prefix
of $h(d)$ for any $d \in \Sigma_4$.

\item[(ii)]
Suppose there exist letters $a, b, c \in \Sigma_4$ and
words $s, t, u, v \in \Sigma_2^*$ such that $h(a) = st$, $h(b) = uv$,
and $h(c) = sv$.  Then either $a = c$ or $b = c$.
\end{itemize}
\label{ming}
\end{fact}

\begin{proof}
\begin{itemize}
\item[(i)]
This can be verified with a short computation.  

\item[(ii)]  This can also be verified with a short computation.
If $|s| \geq 11$, then no two images of distinct
letters share a prefix of length $11$.
If $|s|\leq 10$, then $|t| \geq 11$, and no two images of distinct letters 
share a suffix of length $11$.
\end{itemize}
\endpf
\end{proof}

     Now we resume the proof of Lemma~\ref{seven}.
     For $i = 1, 2, \ldots, n$ define $A_i = h(a_i)$.  
Then if $h(w) = xyyz$, we can write
$$h(w) = A_1 A_2 \cdots A_n = A'_1 A''_1 A_2 \cdots A_{j-1} 
A'_j A''_j A_{j+1} \cdots A_{n-1} A'_n A''_n$$ 
where
\begin{eqnarray*}
A_1 &=& A'_1 A''_1 \\
A_j &=& A'_j A''_j \\
A_n &=& A'_n A''_n \\
x &=& A'_1 \\
y &=& A''_1 A_2 \cdots A_{j-1} A'_j = A''_j A_{j+1} \cdots A_{n-1} A'_n \\
z &=& A''_n, \\
\end{eqnarray*}
and $|A''_1|, |A''_j| > 0$.
See Figure~\ref{fig1}.

\begin{figure}[H]
\begin{center}
\input cube1.pstex_t
\end{center}
\caption{The word $xyyz$ within $h(w)$ \protect\label{fig1}}
\end{figure}

     If $|A''_1| > |A''_j|$, then $A_{j+1} = h(a_{j+1})$ is a factor
of $A''_1 A_2$, hence a factor of $A_1 A_2 = h(a_1 a_2)$.  Thus
we can write $A_{j+2} = A'_{j+2} A''_{j+2}$ with
$$ A''_1 A_2 = A''_j A_{j+1} A'_{j+2}.$$ 
See Figure~\ref{fig2}.

\begin{figure}[H]
\begin{center}
\input cube2.pstex_t
\end{center}
\caption{The case $|A''_1| > |A''_j|$ \protect\label{fig2}}
\end{figure}

But then,
by Fact~\ref{ming} (i), either $|A''_j| = 0$,
or $|A'_{j+2}| = 0$ (so $|A''_1| = |A''_j|$), or $A'_{j+2}$ is a not
a prefix of any $h(d)$.  All three conclusions are impossible.

     If $|A''_1| < |A''_j|$, then $A_2 = h(a_2)$ is a factor of
$A''_j A_{j+1}$, hence a factor of $A_j A_{j+1} = h(a_j a_{j+1})$.
Thus we can write $A_3 = A'_3 A''_3$ with
$$ A''_1 A_2 A'_3 = A''_j A_{j+1} .$$  
See Figure~\ref{fig3}.

\begin{figure}[H]
\begin{center}
\input cube3.pstex_t
\end{center}
\caption{The case $|A''_1| < |A''_j|$ \protect\label{fig3}}
\end{figure}

By Fact~\ref{ming} (i), either $|A''_1| = 0$ or $|A'_j| = 0$ (so $|A''_1| = |A''_j|$)
or $A'_3$ is not a prefix of any $h(d)$.  Again, all three conclusions
are impossible.

     Therefore $|A''_1| = |A''_j|$.  
Hence $A''_1 = A''_j$, $A_2 = A_{j+1}$, $\ldots$, $A_{j-1} = A_{n-1}$,
and $A'_j = A'_n$.  Since $h$ is injective, we have
$a_2 = a_{j+1}, \ldots, a_{j-1} = a_{n-1}$.
It also follows that $|y|$ is divisible by $21$ and
$A_j = A'_j A''_j = A'_n A''_1$.   But by Fact~\ref{ming} (ii), either
(1) $a_j = a_n$ or (2) $a_j = a_1$.  In the first case,
$a_2 \cdots a_{j-1} a_j = a_{j+1} \cdots a_{n-1} a_n$, so
$w$ contains the square $(a_2 \cdots a_{j-1} a_j)^2$, a contradiction.  In the
second case, $a_1 \cdots a_{j-1} = a_j a_{j+1} \cdots a_{n-1}$, so
$w$ contains the square $(a_1 \cdots a_{j-1})^2$, a contradiction.

     This completes the proof of part (i).

    It now remains to prove (ii).  If $h(w)$ contains a ${7 \over 3}^+$-power
$yyy'$, then it contains a square, and by part (i) we know
that $|y| \leq 13$.  
We may assume that $|y'| \leq \lceil {1 \over 3} \cdot 13 \rceil = 5$, so
$|yyy'| \leq 31$.  Hence we need only
check the image of all $36$ squarefree words in $\Sigma_4^3$ to ensure they do not
contain any ${7 \over 3}^+$-power.  We leave this computation to the reader.

\endpf
\end{proof}

      Next, define the substitution $g:\Sigma_3^* \rightarrow 2^{\Sigma_4^*}$
as follows:

\begin{eqnarray*}
g(0) &=& \lbrace 0, 3 \rbrace \\
g(1) &=& \lbrace 1 \rbrace \\
g(2) &=& \lbrace 2 \rbrace .
\end{eqnarray*}

     Let $|w|_a$ denote the number of occurrences of the letter $a$ in the word
$w$.  We prove

\begin{lemma}
Let $w \in \Sigma_3^*$ be any squarefree word.  Then $h(g(w))$ 
is a language of $2^r$ words over $\Sigma_2$, where $r = |w|_0$, and moreover
these words are of length $21 |w|$ and avoid ${7 \over 3}^+$-powers.
\end{lemma}

\begin{proof}
      Let $w$ be a squarefree word over $\Sigma_3$.  Then $g(w)$ is
a language over $\Sigma_4$, and we claim that each word $x \in g(w)$ is squarefree.
For suppose $x \in g(w)$ and $x$ contains a square, say
$x = tuuv$ for some words $t, v \in \Sigma_4^*$, and $u \in \Sigma_4^+$.
Define the morphism $f$ where $f(0) = f(3) = 0$,  $f(1) = 1$, and
$f(2) = 2$.  Then $f(x) = w$ and $f(x)$ contains the square $f(u)f(u)$, a contradiction.
It now follows from Lemma~\ref{seven} that $x$ avoids ${7 \over 3}^+$-powers.
\endpf
\end{proof}

     Finally, we obtain

\begin{theorem}
       Let $C_n$ be the number of binary words of length $n$ that are
${7 \over 3}^+$-power-free.  Then
$C_n = \Omega(\gamma^n)$, where $\gamma = 2^{1/63} \doteq 1.011$.
\end{theorem}

\begin{proof}
      Take any squarefree word $x$ of length $m$ over $\Sigma_3$.  
There must exist a symbol $a \in \Sigma_3$ such that $a$ occurs
at least $\lceil m/3 \rceil$ times in $x$.  By replacing each symbol
$b$ in $x$ with $b-a$ (mod $3$), we get a squarefree word $x'$
with at least $\lceil m/3 \rceil$ occurrences of $0$.

     Now consider $h(g(x'))$.  We get at least
$2^{m/3}$ words of length $21m$, and 
each word is ${7 \over 3}^+$-power-free.  
Write $n = 21m - k$, where $0 \leq k < 21$.   By what precedes,
there are at least $2^{n/63}$ words of length $21m$ that are
${7 \over 3}^+$-power-free.  Thus there are at least
$2^{-k} 2^{n/63} \geq 2^{-21} 2^{n/63}$ words of length $n$ 
with the desired property.
\endpf
\end{proof}

     We have not tried to optimize the value of $\gamma$.  It can be improved
slightly in several ways:  for example, by starting with a squarefree word
over $\Sigma_3$ with a higher proportion of $0$'s; see 
\cite{Tarannikov:2002}.

     For an upper bound on $C_n$, we may reason as follows:  if $w$ is a word
avoiding ${7 \over 3}^+$-powers, then $w$ certainly has no occurrences
of either $000$ or $111$.  Let $E_n$ denote the number of binary words
of length $n$ avoiding both $000$ and $111$.  Then $C_n \leq E_n$ and
it is easy to see that 
\begin{equation}
E_n = E_{n-1} + E_{n-2} 
\label{rec}
\end{equation}
for $n \geq 3$.  Now the characteristic
polynomial of the recursion (\ref{rec}) is $x^2 - x - 1$, and the dominant
zero of this polynomial is $(1 + \sqrt{5})/2 \doteq 1.62$.  By well-known
properties of linear recurrences we get 
$E_n = O(1.62^n)$.

     This procedure may be automated.  Noonan and Zeilberger \cite{Noonan&Zeilberger:1999}
have written a Maple package {\tt DAVID\_IAN} that allows one to specify a list $L$
of forbidden words, and computes the generating function enumerating words avoiding
members of $L$.  We used this package for a list $L$ of $58$ words of length $\leq 24$:
$$ 000,111,01010,10101, \ldots, 110110010011011001001101 $$
including words of the form $x^{(1 + \lfloor 7|x|/3 \rfloor)/|x|}$ for
$1 \leq |x| \leq 10$.  (Words for which shorter members of $L$ are factors can be
omitted.)  We obtained a characteristic polynomial of degree $39$ with
dominant root $1.22990049\cdots$.  Therefore we have shown

\begin{theorem}
The number $C_n$ of binary words of length $n$ avoiding ${7 \over 3}^+$-powers 
satisfies $C_n = O(1.23^n)$.
\end{theorem}

\section{Avoiding arbitrarily large squares}

     Dekking \cite{Dekking:1976} proved that every infinite overlap-free
binary word must contain arbitrarily large squares.  He also
proved that there exists an infinite cubefree binary word
that avoids squares $xx$ with $|x| \geq 4$.  Furthermore the number $4$ is
best possible, since every binary word of length $\geq 30$ contains a cube
or a square $xx$ with $|x| \geq 3$.

     This leads to the following natural question:  what is the largest
exponent $\alpha$ such that every infinite $\alpha$-power-free binary word
contains arbitrarily large squares?  From Dekking's results we know
$2 < \alpha < 3$.    The answer is given in the following theorem.

\begin{theorem}
\begin{itemize}
\item[(i)] Every infinite ${7 \over 3}$-power-free binary word
contains arbitrarily large squares.

\item[(ii)] There exists
an infinite ${7 \over 3}^+$-power-free binary word such that
each square factor $xx$ satisfies $|x| \leq 13$.
\end{itemize}
\label{sq}
\end{theorem}

\begin{proof}
     For (i), let $\bf w$ be an infinite $7 \over 3$-power-free binary
word.  By Theorem~\ref{rs} and Eq.\ (\ref{ineq1}),
any prefix of $\bf w$ of length $2^{n+5}$ contains
$\mu^{n+2} (0)$ as a factor.
But $\mu^{n+2}(0) = \mu^{n}(0110)$,
so any prefix of length $2^{n+5}$ contains the square factor
$xx$ with $x = \mu^n (1)$.

     For (ii), from Theorem~\ref{seven} it follows that
if $\bf w$ is an infinite squarefree word over $\Sigma_4$,
and $h$ is the morphism
defined in \S~\ref{expo-sec}, then $h({\bf w})$ has the desired
properties.
\endpf
\end{proof}

     We note that the number $13$ in Theorem~\ref{sq} (ii)
is not best possible.   A forthcoming paper examines this question
in more detail.

\section{Numerical Results}

     Let $A_n$ (resp., $B_n$, $C_n$, $D_n$) denote the number
of overlap-free words (resp., $7 \over 3$-power-free words,
${7 \over 3}^+$-power-free words, cubefree words)
over the alphabet $\Sigma_2$.  We give here the values of these
sequences for $0 \leq n \leq 28$.

\begin{center}
\begin{tabular}{|r|rrrrrrrrrrrrrrrrr|}
\hline
$n$  &  0 & 1 & 2 & 3 &  4 &  5 &  6 &  7 &  8 &  9 & 10 & 11 &  12 &  13 &  14 & 15 & 16 \\
\hline
$A_n$&  1 & 2 & 4 & 6 & 10 & 14 & 20 & 24 & 30 & 36 & 44 & 48 &  60 &  60 &  62 & 72 & 82 \\
\hline
$B_n$&  1 & 2 & 4 & 6 & 10 & 14 & 20 & 24 & 30 & 40 & 48 & 56 &  64 &  76 &  82 & 92 &106\\
\hline
$C_n$&  1 & 2 & 4 & 6 & 10 & 14 & 20 & 30 & 38 & 50 & 64 & 86 & 108 & 136 & 178 &222 &276\\
\hline
$D_n$&  1 & 2 & 4 & 6 & 10 & 16 & 24 & 36 & 56 & 80 &118 &174 & 254 & 378 & 554 &802&1168\\
\hline
\end{tabular}
\end{center}

\begin{center}
\begin{tabular}{|r|rrrrrrrrrrrr|}
\hline
$n$  & 17 & 18 & 19 & 20 &  21 & 22 & 23 & 24 & 25 & 26 & 27 & 28 \\
\hline
$A_n$& 88 & 96 & 112 & 120 & 120 & 136 & 148 & 164 & 152 & 154 & 148 & 162 \\ 
\hline
$B_n$& 124 & 142 & 152 & 172 & 192 & 210 & 220 & 234 & 256 & 284 & 308 & 314 \\ 
\hline
$C_n$& 330 & 408 & 500 & 618 & 774 & 962 & 1178 & 1432 & 1754 & 2160 & 2660 & 3292 \\ 
\hline
$D_n$& 1716 & 2502 & 3650 & 5324 & 7754 & 11320 & 16502 & 24054 & 35058 & 51144 & 74540 & 108664 \\ 
\hline
\end{tabular}
\end{center}


\section{Acknowledgments}

We would like to thank Ming-wei Wang, who showed the first author
the usefulness of Fact~\ref{ming} (i).

\newcommand{\noopsort}[1]{} \newcommand{\singleletter}[1]{#1}

\end{document}